\documentclass[12pt]{article}
\usepackage{amssymb,amsmath}
\usepackage{cases}
\usepackage{amsfonts}
\usepackage{color,xcolor}
\usepackage[left=2.0cm,right=2.0cm,top=2.0cm,bottom=2.0cm]{geometry}
\usepackage[colorlinks,citecolor=blue,urlcolor=blue]{hyperref}

\newtheorem{theorem}{Theorem}[section]

\newtheorem{lemma}{Lemma}[section]

\newtheorem{definition}{Definition}[section]
\newtheorem{remark}{Remark}[section]

\newcommand{\bal}{\begin{align}}
\newcommand{\bbal}{\begin{align*}}
\newcommand{\beq}{\begin{equation}}
\newcommand{\eeq}{\end{equation}}
\newcommand{\bca}{\begin{cases}}
\newcommand{\eca}{\end{cases}}
\def\div{\mathord{{\rm div}}}
\newcommand{\pa}{\partial}
\newcommand{\fr}{\frac}
\newcommand{\na}{\nabla}
\newcommand{\De}{\dot{\Delta}}

\newcommand{\la}{\lambda}
\newcommand{\cd}{\cdot}
\newcommand{\ep}{\varepsilon}
\newcommand{\dd}{\mathrm{d}}
\newcommand{\B}{\dot{B}}
\newcommand{\LL}{\tilde{L}}
\newcommand{\R}{\mathbb{R}}
\newcommand{\D}{\mathrm{div} }

\newcommand{\uu}{\mathbf{u}}
\newcommand{\vv}{\mathbf{v}}

\begin{document}
\title{Global small solution and optimal decay rate for the Korteweg system in Besov spaces}

\author{Jinlu $\mbox{Li}^{1}$\footnote{E-mail: lijl29@mail2.sysu.edu.cn} \quad \quad \quad Yanghai $\mbox{Yu}^{2}$\footnote{E-mail: yuyanghai214@sina.com}\quad and \quad Weipeng $\mbox{Zhu}^{3}$\footnote{E-mail: mathzwp2010@163.com}\\
\small $^1\mbox{School}$  of Mathematics and Computer Sciences, Gannan Normal University, Ganzhou 341000, China\\
\small $^2\mbox{School}$ of Mathematics and Statistics, Anhui Normal University, Wuhu, Anhui, 241002, China\\
\small $^3\mbox{Department}$ of Mathematics, Sun Yat-sen University, Guangzhou 510275, China}
\date{\today}

\maketitle\noindent{\hrulefill}

{\bf Abstract:} In this paper we consider the Cauchy problem to the Korteweg system with the general pressure in dimension $d\geq2$, and establish the global well-posedness of strong solution for the small initial data in $L^p$ type critical Besov spaces by using the Friedrich method and compactness arguments. Furthermore, we also obtain the optimal decay rate for the Korteweg system in $L^2$ type Besov spaces.

{\bf Keywords:} Korteweg system; global small solution; Besov spaces; optimal decay rate

{\bf MSC (2010):} 35Q35; 76N10
\vskip0mm\noindent{\hrulefill}

\section{Introduction}\label{sec1}
This paper focuses on the Cauchy problem to the following isothermal model of capillary fluids derived by Dunn and Serrin \cite{Dunn 1985}, which can be used as a phase transition model
\begin{equation}\label{eq1.1}\begin{cases}
\ \partial_t \rho+\mathrm{div}(\rho\mathbf{u})=0,& x\in \R^d,t>0, \\
\ \partial_t(\rho\mathbf{u})+\mathrm{div}(\rho\mathbf{u}\otimes\mathbf{u})-\mathrm{div}(\mu(\rho) \mathbb{D}\mathbf{u})-\nabla(\lambda(\rho)\D \uu)+\na P(\rho)=\D \mathbb{K},& x\in \R^d,t>0,\\
\ (\mathbf{u},\rho)|_{t=0}=(\mathbf{u}_0,\rho_0),& x\in \R^d,t>0,
\end{cases}\end{equation}
here $\mathbf{u}(t,x)$ denotes the velocity field and $\rho=\rho(t,x)\in\R^+$ is the density, respectively. $\mu(\rho)$ and $\lambda(\rho)$ are the shear and bulk viscosity coefficients of the flow which fulfill the standard strong parabolicity assumption:
$$\mu>0\quad\mbox{and}\quad2\mu+\lambda>0.$$
In the physical case the viscosity coefficients satisfy $2\mu+d\lambda>0$ which is a special case of the previous assumption. The strain tensor $\mathbb{D}\uu =\frac{1}{2}(\nabla \uu+\nabla^t \uu)$ is the symmetric part of the velocity gradient. The pressure $P$ is suitably smooth increasing function of the density $\rho$.
 The Korteweg tensor $\div\mathbb{K}$ allows to describe the variation of density at the interfaces between two phases, generally a mixture liquid-vapor, which can be written as follows:
\begin{eqnarray}\label{eq1.2}
\div\mathbb{K}=\nabla\Big(\rho\kappa(\rho)\Delta\rho+\frac{\kappa(\rho)+\rho\kappa'(\rho)}{2}|\nabla\rho|^2\Big)-\div(\kappa(\rho)\nabla\rho\otimes\nabla\rho),
        \end{eqnarray}
where the regular function $\kappa$ denotes the capillary coefficient.
It is worth mentioning that a typical pressure term $P(\rho)$ is assumed to obey the following polytropic law,
\begin{eqnarray}\label{eq1.3}
P(\rho) = A\rho^{\gamma},
\end{eqnarray}
where $A$ is the entropy constant and $\gamma\geq1$ is called the adiabatic index.

When $\div\mathbb{K} \equiv 0$, the system (\ref{eq1.1}) reduces to the compressible Navier-Stokes equations (CNS). There have been huge literatures on the study of CNS by many physicists and mathematicians due to its physical importance, complexity, rich phenomena and mathematical challenges(see \cite{C.M.Z,Danchin5 2016,Danchin6 2017,Danchin1 2000,Danchin2 2001,Danchin3 2001} and the references therein). To explore the scaling invariance property, Danchin first introduce in his series papers  the ``Critical Besov Spaces" which were inspired by those efforts on the incompressible Navier-Stokes. By critical, we mean that the solution space that we shall consider has the same scaling invariance by time and space dilations as (\ref{eq1.1}) itself. Precisely speaking, we can check that if $(\rho,\uu)$ solves (\ref{eq1.1}), so does $(\rho_\ell,\uu_\ell)$ where:
$$(\rho_\ell,\uu_\ell)(t,x)=(\rho(\ell^2t,\ell x),\ell \uu(\ell^2t,\ell x)),$$
provided that the pressure term $P$ has been changed into $\ell^2P$. This suggests us to choose initial data $(\rho_0,\uu_0)$ in ``critical spaces" whose norm is invariant for all $\ell>0$ by the transformation $(\rho_0,\uu_0)=(\rho(\ell \cd),\ell \uu(\ell \cd))$. Now, we briefly review some results concerned with the multi-dimensional compressible Korteweg system in the framework of critical Besov spaces which are more relatively with our problem. Danchin and Desjardins in \cite{Danchin 2001} have been the first to obtain the existence of global strong solution with small initial data when $(\rho_0-\bar{\rho}, \uu_0)\in\B^{\fr d2}_{2,1}\times\B^{\fr d2-1}_{2,1}$. Recently, Haspot consider the cases when the viscosity coefficients $\mu(\rho),\la(\rho)$ and the pressure $P(\rho)$ depends linearly on the density for the system (\ref{eq1.1})-(\ref{eq1.2}) with $\kappa(\rho)=\frac{\kappa}{\rho}$, and obtain the global solution with the suitable small initial data in the $L^2$ framework. Subsequently, Haspot continue to investigate the Cauchy problem of the system (\ref{eq1.1})-(\ref{eq1.2}) with $\mu(\rho)=\mu{\rho},\la(\rho)=0$ and the pressure $P(\rho)=\rho$, and establish the global solution under the setting of slightly subcritical $L^p$ type initial data, where the specific choice of the pressure is crucial since it provides a gain of integrability on the effective velocity $\vv$. However, for the general pressure which fulfills (\ref{eq1.3}), the existence of global strong solutions with large initial data is still an open problem even in dimension $d=2$. Motivated by this work, our goal in the present paper is devoted to study the global well-posedness of strong solution for the system (\ref{eq1.1})-(\ref{eq1.2}) with the small initial data in the $L^p$ framework. Concerning the existence of global weak solutions, we refer to \cite{Bresch 2003,Haspot 2011,Jngel 2010} and the references therein.

So far there are lots of mathematical results on the large time behavior of the solution of (CNS) for the initial data $(\rho_0,\uu_0)$ which is a small perturbation of an equilibrium state $(\bar{\rho},0)$. In this research direction, the first effort is due to
Matsumura and Nishida \cite{Matsumura,Matsumura1980}, they showed the global in time existence of solution to (CNS) in $\R^3$ assumed that the initial small perturbation $(\rho_0-\bar{\rho}, \uu_0)$ in $H^3(\R^3) \cap L^1(\R^3)$ and obtained the following decay estimates
$$||\nabla^k(\rho-\bar{\rho},\uu)(t)||_{L^2}\lesssim(1+t)^{-\frac 34-\fr k2}\quad\mbox{for}\quad k=0,1.$$
Subsequently, Ponce \cite{Ponce} obtained more general $L^p$ decay rates
$$||\nabla^k(\rho-\bar{\rho},\uu)(t)||_{L^p}\lesssim(1+t)^{-\frac d2(1-\fr 1p)-\fr k2},\quad2\leq p\leq\infty\quad 0\leq k\leq2 \quad d=2,3.$$
To the best of our knowledge, the approach relies heavily on the analysis of the linearization of the system and the decay estimates of the semi-group which is the reason why almost all the convergence rate of solution as $t\to\infty$ are restricted to the regime which is near the equilibrium state. We also mention those works by \cite{Kawashita,Li,Wang} and the references therein. Okita \cite{Okita} showed the global in time existence of solution to (CNS) for $d\geq3$ when the initial perturbation $(\rho_0-\bar{\rho}, \uu_0)$ is sufficiently small in $((B^{\fr d2}_{2,1}\cap\B_{1,\infty}^0\times(B^{\fr d2-1}_{2,1})\cap\B_{1,\infty}^0))$. Furthermore, Okita apply the decay estimates for the low frequency part of $E(t)$ which is generated by the linearized operator at the constant state $(\bar{\rho}, 0)$ and invoke the energy method for the high frequency part, and obtain the optimal $L^2$ decay rate for strong solutions in critical Besov spaces, namely,
$$||(\rho-\bar{\rho},\uu)(t)||_{B^{\frac d2-1}_{2,1}}\lesssim(1+t)^{-\frac d4}\quad\mbox{for}\quad t\geq0.$$
Recently, Danchin established the optimal time-decay estimate of solution to (CNS) in the general critical Besov framework (for more details, see\cite{Danchin7 2017}). All the decay results we mentioned above are obtained for the system (CNS) which is corresponds to (\ref{eq1.1}) with $\div\mathbb{K}=0$. Inspired by Okita's work, our second purpose is to consider the longtime behavior of the solution of (CNS) with the Korteweg tensor $\div\mathbb{K}$ in the present paper.
\subsection{Reformulation of the System}\label{subsec1.1}
The main difficulties in the study of the compressible fluid flows when dealing with the vacuum is that the momentum equation loses its parabolic regularizing effect. That is why in the present paper we suppose that the initial data $\rho_0$ is a small perturbation of an equilibrium state $\bar{\rho}=1$ (just for convenience). Following the specific choice on the capillary coefficient $\kappa(\rho)=\frac{\kappa}{\rho}$ with $\kappa=\mu^2$ in \cite{Haspot 2017}, then from (\ref{eq1.2}), we have $\div\mathbb{K}=\mu^2\div(\rho\nabla\nabla\ln\rho)$. Denoting $a=\rho-1$ and introducing the effective velocity $\mathbf{v}=\uu+\mu\nabla \ln(1+a)$, hence, as long as $\rho$ does not vanish, we can reformulate the system (\ref{eq1.1})-(\ref{eq1.2}) equivalently as follows
\begin{equation}\label{eq1.4}
\begin{cases}
\ \partial_ta-\mu \Delta a+\mathrm{div}\mathbf{v}=f(a,\vv):=-\mathrm{div}(a\mathbf{v}),  \\
\ \partial_t\mathbf{v}-\mu\Delta \mathbf{v}+\nabla a=\mathbf{g}(a,\vv):=2\mu\nabla(\ln(1+a))\nabla \mathbf{v}-I(a)\na a-\mathbf{v}\cdot\nabla\mathbf{v},\\
\ (a,\mathbf{v})|_{t=0}=(a_0,\mathbf{v}_0),
\end{cases}
\end{equation}
where, we denote $I(a):=\fr{P'(1+a)}{1+a}-1$ and assume that $P'(1)=1$ without loss of generality. Here, denoting $G(a):=\int^a_0(\fr{P'(1+t)}{1+t}-1)\dd t$, we also have $\nabla G(a)=I(a)\nabla a$.

\subsection{Main Results}\label{subsec1.2}
We denote
\begin{align*}
\mathcal{E}=\Big\{(a,\vv)\in \Big(&L^1(\R^+;\mathfrak{\dot{B}}^{\frac d2+1,\frac dp+2}_{2,p})\cap \mathcal{C}(\R^+;\mathfrak{\dot{B}}^{\frac d2-1,\frac dp}_{2,p})\Big)
\\&\times \Big(L^1(\R^+;\mathfrak{\dot{B}}^{\frac d2+1,\frac dp+1}_{2,p})\cap \mathcal{C}(\R^+;\mathfrak{\dot{B}}^{\frac d2-1,\frac dp-1}_{2,p})\Big)^d\Big\}.
\end{align*}
Our first result of this paper as follows.
\begin{theorem}\label{th1.1}
Let $d\geq 2,2\leq p<d$ and $p\leq\min\{4,\frac{2d}{d-2}\}$. Assume that $(\na a_0,\vv_0)\in \B^{\frac dp-1}_{p,1}$ and $(a^\ell_0,\vv^\ell_0)\in \B^{\frac d2-1}_{2,1}$. There exists a constant $\ep$ such that if
$$||(a_0,\vv_0)||^\ell_{\B^{\frac d2-1}_{2,1}}+||(\na a_0,\vv_0)||^h_{\B^{\frac dp-1}_{p,1}}\leq \ep,$$
then \eqref{eq1.4}  has a unique global-in-time solution $(a,\vv)\in \mathcal{E}$. Moreover, there holds for all $T\geq 0$,
$$X_{p,T}(a,\vv)\leq C_0X_{p,0},$$
where
$$X_{p,0}=||(a_0,\vv_0)||^\ell_{\B^{\frac d2-1}_{2,1}}+||(\na a_0,\vv_0)||^h_{\B^{\frac dp-1}_{p,1}}$$
and
$$X_{p,t}(a,\vv)=||(a,\vv)||^\ell_{\LL^\infty_t(\B^{\frac d2-1}_{2,1})\cap L^1_t(\B^{\frac d2+1}_{2,1})}
+||(\na a,\vv)||^h_{\LL^\infty_t(\B^{\frac dp-1}_{p,1})\cap L^1_t(\B^{\frac dp+1}_{p,1})}.$$
\end{theorem}

\begin{remark}\label{re1.1}
When $p=2$ in Theorem \ref{th1.1}, we then have $a_0\in \B^{\frac d2-1}_{2,1}\cap \B^{\frac d2}_{2,1}$ and $\vv_0\in \B^{\frac d2-1}_{2,1}$. Let $d\geq 3$. Assume that $a_0\in B^{\frac d2}_{2,1}$ and $\vv_0\in B^{\frac d2-1}_{2,1}$ satisfying $||a_0||_{B^{\frac d2}_{2,1}}+||\vv_0||_{B^{\frac d2 -1}_{2,1}}\leq \ep$ for small enough $\ep>0$, we also have $a\in \mathcal{C}([0,T];B^{\frac d2}_{2,1}),\ \vv\in \mathcal{C}([0,T];B^{\frac d2-1}_{2,1})$ for all $T\geq 0$.
\end{remark}

We now state our second result of this paper which gives the optimal $L^2$ decay rate for strong solutions.

\begin{theorem}\label{th1.2}
Let $d\geq 3$. Assume that $a_0\in B^{\frac d2}_{2,1}\cap \B^{0}_{1,\infty}$ and $\vv_0\in B^{\frac d2-1}_{2,1}\cap \B^{0}_{1,\infty}$.  There exists a constant $\ep$ such that if
$$M_0:=||a_0||_{B^{\frac d2}_{2,1}\cap \B^{0}_{1,\infty}}+||\vv_0||_{B^{\frac d2-1}_{2,1}\cap \B^{0}_{1,\infty}}\leq \ep , $$
then \eqref{eq1.4} has a unique global-in-time solution $(a,\vv)$. Furthermore, we have for all $t\geq 0$,
$$||a(t)||_{B^{\frac d2}_{2,1}}+||\vv(t)||_{B^{\frac d2-1}_{2,1}}\leq \widetilde{C}_0(1+t)^{-\frac d4}.$$
\end{theorem}

\begin{remark}
By $B^{\frac d2-1}_{2,1}\subset L^2$, the convergence rate of Theorem \ref{th1.2} is optimal.
\end{remark}

The rest of this paper is structured as follows. In Section 2 we present some notions and basic tools. In Section 3 we establish the a priori estimates for the linearized equation of the system (\ref{eq1.1}) which will be crucial in the proof of Theorem \ref{th1.1}. In Section 4 we prove the Theorem \ref{th1.1} by utilizing the Friedrich method and compactness arguments. In Section 5 we complete the proof of Theorem \ref{th1.2} by spectral analysis to the linearized system.

\section{Preliminaries}\label{sec2}
\setcounter{equation}{0}
In this section, we introduce some notations and conventions, and recall some standard theories of Besov space which will be used throughout this paper.
\subsection{Notations}\label{sec2.1}
In the following, we denote by $(\cdot|\cdot)$ the $L^2$ scalar product and use the convention that $C$, with or without subscripts, to denote strictly positive constants whose values are insignificant and may change from line to line. $A\lesssim B$ means that there is a uniform positive constant $c$ independent of $A$ and $B$ such that $A\leq cB$. For $X$ a Banach space and $T>0$, we denote by $L_T^p(X)$ with $p\in[1,\infty]$ stands for the set of measurable functions on $[0,T]$ with values in $X$, such that $t\mapsto||f(t)||_{X}\in L^p([0,T])$.

\subsection{Littlewood-Paley theory and Besov spaces}\label{sec2.2}
The following material involving the theories of Littlewood-Paley is standard, we refer the readers to Bahouri, Chemin and Danchin \cite{Bahouri2011}.

Let $\mathcal{C}$ denote the annulus $\{\xi\in \R^d:3/4\leq|\xi|\leq8/3\}$ and $\mathcal{B}$ denote the ball $\{\xi\in \R^d:|\xi|\leq4/3\}$. There exist two radial functions $\chi\in C_c^{\infty}(\mathcal{B}(0,4/3))$ and $\varphi\in C_c^{\infty}(\mathcal{C})$ both taking values in $[0,1]$ such that
\begin{eqnarray*}
\sum_{j\in\mathbb{Z}}\varphi(2^{-j}\xi)=1 \quad \mbox{for} \quad\xi\in \R^d\setminus\{0\} \quad \mbox{and} \quad \chi(\xi)+\sum_{j\geq0}\varphi(2^{-j}\xi)=1 \quad \mbox{for} \quad \xi\in \R^d.
\end{eqnarray*}
For every $f\in \mathcal{S'}(\mathbb{R}^d)$, the homogeneous (or nonhomogeneous) dyadic blocks ${\dot{\Delta}}_j$ (or $\Delta_j$) and homogeneous (or nonhomogeneous) low-frequency cut-off operator $\dot{S}_j$ (or $S_j$) are defined as follows
$$
\forall j\in \mathbb{Z}, \; {\dot{\Delta}}_jf=\varphi(2^{-j}D)f \quad \mbox{and}\quad\dot{S}_jf=\chi(2^{-j}D)f=\sum_{q\leq j-1}{\dot{\Delta}}_qf;
$$
\begin{numcases}{\Delta_jf=}
0, & $j\leq-2$;\nonumber\\
\chi(D)f, & $j=-1$;\nonumber\\
\varphi(2^{-j}D)f, & $j\geq0$;\nonumber
\end{numcases}
and
$$
 S_jf=\sum_{q=-1}^{j-1}{\Delta}_qf.
$$
Unfortunately, for the homogeneous case, the Littlewood-Paley decomposition is invalid. We need to a new space to modify it, namely,
\begin{eqnarray*}
\mathcal{S}'_h\triangleq\Big\{f \in \mathcal{S'}(\mathbb{R}^{d}):\; \lim_{j\rightarrow-\infty}||\chi(2^{-j}D)f||_{L^{\infty}}=0 \Big\}.
\end{eqnarray*}
Then we have the formal Littlewood-Paley decomposition in the homogeneous case
\begin{eqnarray*}
f=\sum_{j\in\mathbb{Z}}\dot{\Delta}_jf, \quad \forall f\in \mathcal{S}'_h.
\end{eqnarray*}
With a suitable choice of $\varphi$, one can easily verify that the Littlewood-Paley decomposition satisfies the property of almost orthogonality:
 \begin{eqnarray}\label{lyz-1}
\dot{\Delta}_j\dot{\Delta}_kf\equiv0 \quad \mbox{if} \quad |j-k|\geq2 \quad \mbox{and} \quad \dot{\Delta}_j(\dot{S}_{k-1}f\dot{\Delta}_kf)\equiv0 \quad \mbox{if} \quad |j-k|\geq5.
\end{eqnarray}
Next we recall Bony's decomposition from \cite{Bahouri2011}:
\begin{align*}
uv=\dot{T}_uv+\dot{T}_vu+\dot{R}(u,v),
\end{align*}
with
\begin{align}\label{lyz-2}
\dot{T}_uv=\sum_{j\in\mathbb{Z}}\dot{S}_{j-1}u\dot{\Delta}_jv, \quad \quad \dot{R}(u,v)=\sum_{j\in\mathbb{Z}}\dot{\Delta}_ju\widetilde{\Delta}_jv, \quad \quad \widetilde{\Delta}_jv=\sum_{|j'-j|\leq 1}\dot{\Delta}_{j'}v.
\end{align}

The operators $\dot{\Delta}_j$ and $\Delta_j$ help us recall the definition of the inhomogenous Besov spaces, the homogenous Besov spaces and hybrid-Besov spaces (see \cite{Bahouri2011})

\begin{definition}
Let $s\in \mathbb{R},T>0$ and $1\leq p,r,q\leq \infty$. The homogeneous Besov space $\B^s_{p,r}$ is the set of tempered distribution $f\in \mathcal{S}'_h$ satisfying
\begin{align*}
||f||_{\B^s_{p,r}}\triangleq \Big|\Big|(2^{js}||\dot{\Delta}_j f||_{L^p})_{j}\Big|\Big|_{\ell^r}<\infty.
\end{align*}
The time-sapce homogeneous Besov spaces is the set of tempered distribution $f$ satisfying
\bbal
\lim_{j\rightarrow -\infty}||\dot{S}_jf||_{L^q_T(\LL^\infty)}=0
\end{align*}
and
\begin{align*}
||f||_{\tilde{L}^q_T(\B^s_{p,r})}\triangleq \Big|\Big|\big(2^{js}||\dot{\Delta}_jf(t)||_{L^q_T(L^p)}\big)_j \Big|\Big|_{\ell^r}<+\infty.
\end{align*}
The nonhomogeneous Besov space $B^s_{p,r}$ is the set of tempered distribution $f$ satisfying
\begin{align*}
||f||_{B^s_{p,r}}\triangleq \Big|\Big|(2^{js}||\Delta_j f||_{L^p})_{j}\Big|\Big|_{\ell^r}<\infty.
\end{align*}
\end{definition}

\begin{remark}\label{rem2.1}
\begin{itemize}
  \item Restricting the above norms to the low or high frequencies parts of distributions will be fundamental in our approach. Fix some integer $j_0$, we denote
\bbal
&||f||^\ell_{\B^s_{p,1}}=\sum_{j\leq j_0}2^{js}||\De_jf||_{L^p}\quad  \mbox{and}\quad  ||f||^h_{\B^s_{p,1}}=\sum_{j\geq j_0+1}2^{js}||\De_jf||_{L^p};
\\&||f||^\ell_{\tilde{L}^q_T(\B^s_{p,1})}=\sum_{j\leq j_0}2^{js}||\De_jf||_{L^q_T(L^p)}\quad  \mbox{and}\quad  ||f||^h_{\tilde{L}^q_T(\B^s_{p,1})}=\sum_{j\geq j_0+1}2^{js}||\De_jf||_{L^q_T(L^p)}.
\end{align*}
  \item By Minkowski's inequality, it is easy to find that
\begin{align*}
||f||_{\tilde{L}^q_T(\B^s_{p,r})}\leq ||f||_{L^q_T(\B^s_{p,r})} \quad \mathrm{if} \quad q\leq r\quad \mbox{and} \quad ||f||_{\tilde{L}^q_T(\B^s_{p,r})}\geq ||f||_{L^q_T(\B^s_{p,r})} \quad \mathrm{if} \quad q\geq r.
\end{align*}
\end{itemize}
\end{remark}
Now we introduce the hybrid-Besov space we will work with in this paper. Let $j_0>0$ be as in Lemma \ref{lem3.1}.
\begin{definition}
Let $s,\sigma\in \mathbb{R}$ and $1\leq p\leq \infty$. The hybrid-Besov space $\mathfrak{\dot{B}}^{s,\sigma}_{2,p}$ is the set of tempered distribution $f\in \mathcal{S}'_h$ satisfying
\begin{align*}
||f||_{\mathfrak{\dot{B}}^{s,\sigma}_{2,p}}\triangleq \sum_{j\leq j_0}2^{js}||\De_jf||_{L^2}+\sum_{j\geq j_0+1}2^{j\sigma}||\De_jf||_{L^p}<\infty.
\end{align*}
The time-space hybrid-Besov space $\LL^q_T(\mathfrak{\dot{B}}^{s,\sigma}_{2,p})$ is the set of tempered distribution $f$ satisfying
\bbal
\lim_{j\rightarrow -\infty}||\dot{S}_jf||_{L^q_T(\LL^\infty)}=0,
\end{align*}
and
\begin{align*}
||f||_{\LL^q_T(\mathfrak{\dot{B}}^{s,\sigma}_{2,p})}\triangleq \sum_{j\leq j_0}2^{js}||\De_jf||_{L^q_T(L^2)}+\sum_{j\geq j_0+1}2^{j\sigma}||\De_jf||_{L^q_T(L^p)}<\infty.
\end{align*}
\end{definition}

\subsection{Basic Properties}\label{sec2.3}
The following Bernstein lemma will be stated as follows (see \cite{Bahouri2011}):
\begin{lemma}(\cite{Bahouri2011})\label{lem2.1}
Let $1\leq p\leq q\leq \infty$  and $\mathcal{B}$ be a ball and $\mathcal{C}$ a ring of $\mathbb{R}^d$. Assume that $f\in L^p$, then for any $\alpha\in\mathbb{N}^d$, there exists a constant $C$ independent of $f$, $j$ such that
\begin{align*}
&\mathrm{Supp} \,\hat{f}\subset\lambda \mathcal{B}\Rightarrow\sup_{|\alpha|=k}\|\partial^{\alpha}f\|_{L^q}\le C^{k+1}\lambda^{k+d(\frac1p-\frac1q)}\|f\|_{L^p},\\
&\mathrm{Supp} \,\hat{f}\subset\lambda \mathcal{C}\Rightarrow C^{-k-1}\lambda^k\|f\|_{L^p}\le\sup_{|\alpha|=k}\|\partial^{\alpha}f\|_{L^p}
\le C^{k+1}\lambda^{k}\|f\|_{L^p}.
\end{align*}
\end{lemma}
As a result of Bernstein's inequalities, we have the following Besov embedding theorem.
\begin{lemma}(\cite{Bahouri2011})\label{lem2.2} Let $1\leq p_1\leq p_2\leq\infty$ and $1\leq r_1\leq r_2\leq\infty$. Then we have for $s\in\R$
$$\dot{B}_{p_1,r_1}^{s}\hookrightarrow \dot{B}_{p_2,r_2}^{s-d(\frac{1}{p_1}-\frac{1}{p_2})}.$$
\end{lemma}
Next, we give the important product acts on homogenous Besov spaces and composition estimate which will be also often used implicity throughout the paper.
\begin{lemma}(\cite{C.M.Z})\label{lem2.3}
Let $d\geq 2,1\leq p\leq \infty$ and $s\leq d/p,t\leq d/p$ with $s+t>d\max\{0,\frac2p-1\}$. Then we have for $(f,g)\in \B^{s}_{p,1}(\R^d)\times\B^{t}_{p,1}(\R^d)$
\begin{align*}
||fg||_{\B^{s+t-d/p}_{p,1}}&\leq C||f||_{\B^{s}_{p,1}}||g||_{\B^{t}_{p,1}}.
\end{align*}
\end{lemma}

\begin{lemma}(\cite{Bahouri2011})\label{lem2.4}
Let $s>0$, $p\in[1,\infty]$ and $f,g\in \dot{B}_{p,1}^{s}(\R^d)\cap L^{\infty}(\R^d)$.
\begin{enumerate}
  \item Assume that $F\in W_{loc}^{[s]+2,\infty}(\R^d)$ with $F(0)=0$, then there exists a function $C$ depending only on $s, p, d$ and $F$ such that
\begin{align*}
||F(f)||_{\dot{B}_{p,1}^{s}}\leq C(||f||_{L^{\infty}})||f||_{\dot{B}_{p,1}^{s}}.
\end{align*}
  \item Assume that $H\in W_{loc}^{[s]+3,\infty}(\R^d)$ with $H'(0)=0$ and $g\in \dot{B}_{p,1}^{s}(\R^d)\cap L^{\infty}(\R^d)$, then there exists a function $C$ depending only on $s, p, d$ and $H$ such that
\begin{align*}
||H(f)-H(g)||_{\dot{B}_{p,1}^{s}}&\leq C(||f||_{L^{\infty}},||g||_{L^{\infty}})||f-g||_{\dot{B}_{p,1}^{s}\cap L^\infty}\Big(||f||_{\dot{B}_{p,1}^{s}\cap L^\infty}+||g||_{\dot{B}_{p,1}^{s}\cap L^\infty}\Big).
\end{align*}
\end{enumerate}
\end{lemma}
\begin{remark}\label{re2}
Let $s>0$, $p\in[1,\infty]$ and $f,g\in \dot{B}_{p,1}^{s}(\R^d)\cap L^{\infty}(\R^d)$ with $||f||_{L^\infty}\leq \frac12,||g||_{L^\infty}\leq \frac12$. Let $\delta(x)$ be a smooth function with value in $[0,1]$, supported in the ball $B(0,\frac23)$ and equal to 1 on $B(0,\frac35)$. If $E(x)=\ln(1+x)\delta(x)$, then there exists a function $C$ depending only on $s, p, d$ and $H$ such that
\begin{align*}
&||E(f)||_{\dot{B}_{p,1}^{s}}\leq C(||f||_{L^{\infty}})||f||_{\dot{B}_{p,1}^{s}},
\\&||E(f)-E(g)||_{\dot{B}_{p,1}^{s}}\leq C(||f||_{L^{\infty}},||g||_{L^{\infty}})||f-g||_{\dot{B}_{p,1}^{s}\cap L^\infty}\Big(||f||_{\dot{B}_{p,1}^{s}\cap L^\infty}+||g||_{\dot{B}_{p,1}^{s}\cap L^\infty}\Big).
\end{align*}
\end{remark}

Finally, we end this section with the following interpolation inequality.
\begin{lemma}\label{lem2.5} (\cite{Bahouri2011})
For $(p,r_1,r_2,r)\in[1,\infty]^4$, $s_1\neq s_2$ and $\theta\in(0,1)$, the following inequality holds
$$\|u\|_{\dot{B}_{p,r}^{\theta s_1+(1-\theta)s_2}}\leq C\|u\|^{\theta}_{\dot{B}_{p,r_1}^{s_1}}\|u\|^{1-\theta}_{\dot{B}_{p,r_2}^{s_2}}.$$
\end{lemma}

\subsection{Some useful lemmas}\label{sec2.4}

\begin{lemma}(\cite{Okita})\label{lem2.6}
(i) Let $a,b>0$ satisfying $\max\{a,b\}>1$. Then
$$\int^t_0(1+s)^{-a}(1+t-s)^{-b}\mathrm{d} s<C(1+t)^{-\min\{a,b\}}, \qquad t>0.$$
(ii) Let $a,b>0$ and $f\in L^1(0,\infty)$. Then
$$\int^t_0(1+s)^{-a}(1+t-s)^{-b}f\mathrm{d} s<C(1+t)^{-\min\{a,b\}}\int^t_0|f|\mathrm{d} s, \qquad t>0.$$
\end{lemma}

\begin{lemma}(\cite{Danchin 2010})\label{lem2.7}
If Supp $\mathcal{F}f\subset \{\xi\in\R^d:R_1\lambda\leq |\xi|\leq R_2\lambda\}$, then there exists $c$ depending only on $d$, $R_1$ and $R_2$ so that for all $1<p<\infty$,
$$c\lambda^2\Big(\frac{p-1}{p}\Big)\int_{\R^d}|f|^p\dd x\leq (p-1)\int_{\R^d}|\na f|^2|f|^{p-2}\dd x=-\int_{\R^d}\Delta f|f|^{p-2}\dd x.$$
\end{lemma}

\begin{lemma}\label{lem2.8}
Let $d\geq 2,2\leq p<2d$ and $p\leq\min\{4,\frac{2d}{d-2}\}$. Then for $f,g\in \B^{d/p}_{p,1}(\R^d)\cap\B^{d/p-1}_{p,1}(\R^d)$, we have
\begin{align*}
||fg||_{\B^{d/2-1}_{2,1}}&\leq C\big(||f||_{\B^{d/p}_{p,1}}||g||_{\B^{d/p-1}_{p,1}}+||g||_{\B^{d/p}_{p,1}}||f||_{\B^{d/p-1}_{p,1}}\big).
\end{align*}
\end{lemma}
{\bf Proof}\quad In order to prove our claim, by Bony's decompose, we write $fg$ as follows
$$fg=\dot{T}_fg+\dot{T}_gf+\dot{R}(f,g).$$
Let $\frac1q=\frac12-\frac1p$. Since $p\leq\min\{4,\frac{2d}{d-2}\}$, we have $2\leq p\leq 4\leq q$ and $q\geq d$.

For the term $\dot{T}_fg$, we deduce that
\begin{align*}
||\dot{T}_fg||_{\B^{d/2-1}_{2,1}}&\lesssim \sum_{j\in \mathbb{Z}}\sum_{|j-k|\leq 4}2^{j(\frac d2-1)}||\De_j\big(\dot{S}_{k-1}f\dot{\Delta}_{k}g\big)||_{L^2}
\\&\lesssim \sum_{j\in \mathbb{Z}}\sum_{|j-k|\leq 4}2^{j(\frac d2-1)}||\dot{S}_{k-1}f||_{L^q}||\dot{\Delta}_{k}g||_{L^p} \quad (\mbox{by H\"{o}lder's inequality})
\\&\lesssim\sum_{j\in \mathbb{Z}}\sum_{|j-k|\leq 4}2^{j(\frac d2-1)}\sum_{l\leq k-2}||\dot{\Delta}_lf||_{L^{q}}||\dot{\Delta}_kg||_{L^p}
\\&\lesssim\sum_{j\in \mathbb{Z}}\sum_{|j-k|\leq 4}2^{j(\frac d2-1)}\sum_{l\leq k-2}2^{(\frac{d}{p}+1-\frac{d}{2})l}2^{(\frac{d}{p}-1)l}||\dot{\Delta}_lf||_{L^{p}}||\dot{\Delta}_kg||_{L^p}\quad (\mbox{by Bernstein inequaltiy})
\\&\lesssim ||f||_{\B^{d/p-1}_{p,1}}\sum_{j\in \mathbb{Z}}\sum_{|j-k|\leq 4}2^{(j-k)(\frac d2-1)}2^{k{\frac dp}}||\dot{\Delta}_{k}g||_{L^p}
\\&\lesssim ||f||_{\B^{d/p-1}_{p,1}}||g||_{\B^{d/p}_{p,1}},
\end{align*}
where we have used that $p\leq4$ in the fourth step and $p\leq2d/(d-2)$ in the fifth step.

Similarly, for the term $\dot{T}_gf$, we have
\begin{align*}
||\dot{T}_gf||_{\B^{d/2-1}_{2,1}}&\lesssim ||g||_{\B^{d/p-1}_{p,1}}||f||_{\B^{d/p}_{p,1}}.
\end{align*}
By \eqref{lyz-1}-\eqref{lyz-2}, we see that the Fourier transform of $\dot{\Delta}_{k}u\widetilde{\Delta}_{k}v$ is supported in $2^{k}B(0,8)$, which implies
$$\De_j(\dot{\Delta}_{k}u\widetilde{\Delta}_{k}v)=0 \qquad \mathrm{for} \qquad j\geq k+4.$$
For the last term $\dot{R}(f,g)$, we get
\begin{align*}
||\dot{R}(f,g)||_{\B^{d/2-1}_{2,1}}&\lesssim \sum_{j\in \mathbb{Z}}\sum_{k\geq j-3}2^{j(\frac d2-1)}||\dot{\Delta}_j(\dot{\Delta}_{k}f\widetilde{\dot{\Delta}}_kg)||_{L^2}
\\&\lesssim \sum_{j\in \mathbb{Z}}\sum_{k\geq j-3}2^{j(\frac d2-1)}||\dot{\Delta}_{k}f||_{L^{q}}||\widetilde{\dot{\Delta}}_kg||_{L^p}\\
&\lesssim \sum_{j\in \mathbb{Z}}\sum_{k\geq j-3}2^{j(\frac d2-1)}2^{(\frac{2}{p}-\frac{1}{2})dk}||\dot{\Delta}_kf||_{L^{p}}\sum_{|l-k|\leq 1}2^{-\frac{d}{p}l}2^{\frac{d}{p}l}||\dot{\Delta}_lg||_{L^p}\nonumber\\
&\lesssim\sum_{j\in \mathbb{Z}}\sum_{k\geq j-3}2^{j(\frac d2-1)}2^{(\frac{d}{p}-\frac{d}{2})k}||\dot{\Delta}_kf||_{L^{p}}||g||_{\dot{B}_{p,1}^{d/p}}\nonumber\\
&\lesssim\sum_{j\in \mathbb{Z}}\sum_{k\geq j-3}2^{(\frac{d}{2}-1)(j-k)}2^{(\frac{d}{p}-1)k}||\dot{\Delta}_kf||_{L^{p}}||g||_{\dot{B}_{p,1}^{d/p}}
\\&\lesssim ||f||_{\B^{d/p-1}_{p,1}}||g||_{\B^{d/p}_{p,1}}.
\end{align*}
This completes the proof of this lemma.

\section{A priori estimates for the linearized equation}
\setcounter{equation}{0}
In this section, we consider the following linearized equation of the system  which paly an important role in the proof of our theorem:
\begin{equation}\label{eq3.0}
\begin{cases}
\ \partial_ta-\mu \Delta a+\mathrm{div}\mathbf{v}=f,  \\
\ \partial_t\mathbf{v}-\mu\Delta \mathbf{v}+\nabla a=\mathbf{g},\\
\ (a,\mathbf{v})|_{t=0}=(a_0,\mathbf{v}_0).
\end{cases}
\end{equation}

We have the following lemma.
\begin{lemma}\label{lem3.1}
Let $(a,\mathbf{v})$ be the smooth solution of system \eqref{eq3.0}. Then there  holds for all $t\in[0,T]$,
\begin{align*}
&||(a,\mathbf{v})||^\ell_{\LL^\infty_T(\B^{\frac d2-1}_{2,1})\cap L^1_T(\B^{\frac d2+1}_{2,1})}+||(\na a,\mathbf{v})||^h_{\LL^\infty_T(\B^{\frac dp-1}_{p,1})\cap L^1_T(\B^{\frac dp+1}_{p,1})}
\\&\lesssim ||(a_0,\mathbf{v}_0)||^\ell_{\B^{\frac d2-1}_{2,1}}+||(\nabla a_0,\mathbf{v}_0)||^h_{\B^{\frac dp-1}_{p,1}}+||(f,\mathbf{g})||^\ell_{L^1_T(\B^{\frac d2-1}_{2,1})}+||(\na f,\mathbf{g})||^h_{L^1_T(\B^{\frac dp-1}_{p,1})}.
\end{align*}
\end{lemma}
{\bf Proof}\quad Applying the operator $\De_j$ to Equations $\eqref{eq3.0}_1$ and $\eqref{eq3.0}_2$, respectively, yields that
\begin{equation}\label{eq3.1}
\begin{cases}
\partial_ta_j-\mu \Delta a_j+\mathrm{div}\mathbf{v}_j=f_j,  \\
\partial_t\mathbf{v}_j-\mu\Delta \mathbf{v}_j+\nabla a_j=\mathbf{g}_j,
\end{cases}
\end{equation}
here and in the sequel, we always denote $\phi_j=\De_j\phi$.

Taking the $L^2$ inner product of Equations $\eqref{eq3.1}_1$ and $\eqref{eq3.1}_2$ with $a_j$ and $\vv_j$, respectively, then integrating by parts, we get
\begin{equation}\label{e3.1}
\begin{cases}
\frac12 \frac{\dd}{\dd t}||a_j||^2_{L^2}+\mu ||\nabla a_j||^2_{L^2}+(\D \vv_j|a_j)=(f_j|a_j),  \\
\frac12 \frac{\dd}{\dd t}||\vv_j||^2_{L^2}+\mu ||\nabla \vv_j||^2_{L^2}+(\na a_j|\vv_j)=(\mathbf{g}_j|\vv_j).
\end{cases}
\end{equation}
Noticing the fact that $(\D \vv_j|a_j)=-(\na a_j|\vv_j)$ and using Lemma \ref{lem2.1}, then we get from \eqref{e3.1} that
\begin{equation*}
\frac{\dd}{\dd t}||(a_j,\vv_j)||^2_{L^2}+c_0\mu2^{2j}||(a_j,\vv_j)||^2_{L^2}\leq C||(f_j,\mathbf{g}_j)||_{L^2}||(a_j,\vv_j)||_{L^2},
\end{equation*}
which leads to
\begin{equation}\label{eq3.2}
\frac{\dd}{\dd t}||(a_j,\vv_j)||_{L^2}+c_0\mu2^{2j}||(a_j,\vv_j)||_{L^2}\leq C||(f_j,\mathbf{g}_j)||_{L^2}.
\end{equation}
Multiplying both sides of \eqref{eq3.2} by $2^{j(\frac d2-1)}$ and summing up over $j\leq j_0$, we infer that
\bal\label{con1}
&||(a,\mathbf{v})||^\ell_{\LL^\infty_T(\B^{\frac d2-1}_{2,1})\cap L^1_T(\B^{\frac d2+1}_{2,1})}
\leq C||(a_0,\mathbf{v}_0)||^\ell_{\B^{\frac d2-1}_{2,1}}+C||(f,\mathbf{g})||^\ell_{L^1_T(\B^{\frac d2-1}_{2,1})}.
\end{align}
Now, we need to bound the high frequencies of the solution in $L^p$ framework. First, multiplying each component of Equation $\eqref{eq3.1}_2$ by $\vv^i_j|\vv^i_j|^{p-2}$ and integrating over $\R^d$ gives for $i=1,\cdots, d$,
\begin{equation*}
\frac1p \frac{\dd}{\dd t}||\vv^i_j||^p_{L^p}+\int\mu \Delta \vv^i_j\cdot\vv^i_j|\vv^i_j|^{p-2}\dd x+\int\pa_i a_j\vv^i_j|\vv^i_j|^{p-2}\dd x=\int\mathbf{g}^i_j\vv^i_j|\vv^i_j|^{p-2}\dd x.
\end{equation*}
After summation on $i=1,\cdots d$, we conclude from Lemma \ref{lem2.1}  and Lemma \ref{lem2.7} that
\begin{equation*}
\frac1p \frac{\dd}{\dd t}||\vv_j||^p_{L^p}+c_p\mu 2^{2j}||\vv_j||^p_{L^p}\leq C(||\na a_j||_{L^p}+||\mathbf{g}_j||_{L^p})||\vv_j||^{p-1}_{L^p},
\end{equation*}
which leads to
\begin{equation}\label{e3.3}
\frac{\dd}{\dd t}||\vv_j||_{L^p}+c_p\mu 2^{2j}||\vv_j||_{L^p}\leq C||\mathbf{g}_j||_{L^p}+C||\na a_j||_{L^p}.
\end{equation}
Applying the operator $\pa_i$ to the first equation of \eqref{eq3.1}, we get
\begin{equation}\label{eq3.3}
\partial_t\pa_i a_j-\mu \Delta \pa_ia_j+\pa_i\mathrm{div}\mathbf{v}_j=\pa_if_j.
\end{equation}
Then, multiplying \eqref{eq3.3} by $|\pa_ia_j|^{p-2}\pa_ia_j$ and integrating over $\R^d$ gives for $i=1,\cdots d$
\bbal
&\quad \ \frac1p \frac{\dd}{\dd t}||\pa_ia_j||^p_{L^p}+\mu \int\Delta \pa_ia_j\pa_ia_j|\pa_ia_j|^{p-2}\dd x+\int\pa_i \D \vv\pa_ia_j|\pa_ia_j|^{p-2}\dd x
=\int\pa_if_j\pa_ia_j|\pa_ia_j|^{p-2}\dd x.
\end{align*}
After summation on $i=1,\cdots d$, we infer from Lemma \ref{lem2.1} and Lemma \ref{lem2.7} that
\bbal
\frac1p \frac{\dd}{\dd t}||\na a_j||^p_{L^p}+c_p\mu 2^{2j}||\na a_j||^{p}_{L^p}\leq C(||\na f_j||_{L^p}+2^{2j}||\vv_j||_{L^p})||\na a_j||^{p-1}_{L^p}.
\end{align*}
which implies
\bal\label{e3.4}
 \frac{\dd}{\dd t}||\na a_j||_{L^p}+c_p\mu 2^{2j}||\na a_j||_{L^p}\leq C||\na f_j||_{L^p}+C2^{2j}||\vv_j||_{L^p}.
\end{align}
Hence, adding the inequality $\eqref{e3.4}\times\gamma$ to the inequality \eqref{e3.3}, we find that for $j\geq j_0+1$,
\begin{align}\label{eq3.4}
&\frac{\dd}{\dd t}(\gamma||\na a_j||_{L^p}+||\vv_j||_{L^p})+c_p\mu2^{2j}||\vv_j||_{L^p}+c_p\mu\gamma2^{2j}||\na a_j||_{L^p}\nonumber\\
\leq& C\gamma 2^{2j}||\vv_j||_{L^p}+C2^{-2j_0}2^{2j}||\na a_j||_{L^p}+C||(\gamma\na f_j,\mathbf{g}_j)||_{L^p}.
\end{align}
Choosing $\gamma$ suitably small and $j_0$ sufficiently large and absorbing the first two terms of RHS of \eqref{eq3.4}, we discover that
\begin{align}\label{e3.5}
&\quad \ \frac{\dd}{\dd t}(\gamma||\na a_j||_{L^p}+||\vv_j||_{L^p})+2^{2j}c_p\mu(||\vv_j||_{L^p}+\gamma||\na a_j||_{L^p})\leq C||(\na f_j,\mathbf{g}_j)||_{L^p}.
\end{align}
Hence, multiplying \eqref{e3.5} by $2^{j(\frac dp-1)}$ and summing up over $j\geq j_0+1$ yields that
\bal\label{con2}
&||(\na a,\mathbf{v})||^h_{\LL^\infty_T(\B^{\frac dp-1}_{p,1})\cap L^1_T(\B^{\frac dp+1}_{p,1})}
\leq C||(\na a_0,\mathbf{v}_0)||^h_{\B^{\frac dp-1}_{p,1}}+C||(\na f,\mathbf{g})||^h_{L^1_T(\B^{\frac dp-1}_{p,1})}.
\end{align}
Therefore, combining \eqref{con1} and \eqref{con2}, we complete the proof of lemma \ref{lem3.1}.

\section{Proof of the main theorem}
\setcounter{equation}{0}
Now, we will divide the proof of Theorem \ref{th1.1} into several steps. The method of the proof is a very classical one.

\textbf{Step 1: Building of the approximation sequence.}

Let $L^2_n$ be the set of $L^2$ function spectrally supported in the Balls $2^n\mathcal{B}$ and let $\Omega_n$ be the set of functions $(a,\mathbf{v})$ of $(L^2_n)^{1+d}$ such that $||a(x)||_{L^\infty(\mathbb{R}^2)}\leq \frac12$. Let us consider the following approximate system
\begin{equation}\label{eq4.1}
\begin{cases}
\ \partial_ta^n-\mu \Delta a^n+\mathrm{div}\mathbf{v}^n=-S_n\big(\mathrm{div}(a^n\mathbf{v}^n)\big):=f^n,  \\
\ \partial_t\mathbf{v}^n-\mu\Delta \mathbf{v}^n+\nabla a^n=S_n\big(2\mu\nabla(\ln(1+a^n))\nabla \mathbf{v}^n-I(a^n)\na a^n-\vv^n\cdot\nabla\vv^n\big):=g^n,\\
\ (a^n,\mathbf{v}^n)|_{t=0}=S_n(a_0,\mathbf{v}_0).
\end{cases}
\end{equation}
It is easy to show that system \eqref{eq4.1} has a unique solution $(a^n,\mathbf{v}^n)_{n\in \mathbb{N}}$ in the space $\mathcal{C}^1([0,T^*_n);\Omega_n)$.

\textbf{Step 2: Uniform estimates. }

Applying Lemma \ref{lem3.1} to system \eqref{eq4.1}, we have for any $T\in[0,T^*_n)$
\bal\label{e4.2}
X_{p,T}(a^n,\vv^n)\lesssim X_{p,0}+||(f^n,\mathbf{g}^n)||^\ell_{L^1_T(\B^{\frac d2-1}_{2,1})}+||(\na f^n,\mathbf{g}^n)||^h_{L^1_T(\B^{\frac dp-1}_{p,1})}.
\end{align}
According to the definition for the Besov spaces, it is easy to show that
\bbal
||(a^n,\vv^n)||_{\LL^\infty_T(\B^{\frac dp-1}_{p,1})\cap L^1_T(\B^{\frac dp+1}_{p,1})}+||a^n||_{\LL^\infty_T(\B^{\frac dp}_{p,1})\cap L^1_T(\B^{\frac dp+2}_{p,1})}\lesssim X_{p,T}(a^n,\vv^n).
\end{align*}
In order to bound the high frequency of $f^n$ and $\mathbf{g}^n$, we see from  Lemmas \ref{lem2.3}-\ref{lem2.4}, Remark \ref{re2} and Lemma \ref{lem2.8} that
\bbal
&\quad ||\D(a^n\vv^n)||^\ell_{L^1_T(\B^{\frac d2-1}_{2,1})}+||\D(a^n\vv^n)||^h_{L^1_T(\B^{\frac dp}_{p,1})}
\\&\lesssim ||\vv^n||_{\LL^\infty_T(\B^{\frac dp-1}_{p,1})}||a^n||_{L^1_T(\B^{\frac dp+1}_{p,1})}+||\vv^n||_{\LL^2_T(\B^{\frac dp}_{p,1})}||a^n||_{\LL^2_T(\B^{\frac dp}_{p,1})}
+||a^n||_{\LL^\infty_T(\B^{\frac dp-1}_{p,1})}||\vv^n||_{L^1_T(\B^{\frac dp+1}_{p,1})}
\\&\quad+||\vv^n||_{\LL^2_T(\B^{\frac dp}_{p,1})}||a^n||_{\LL^2_T(\B^{\frac dp+1}_{p,1})}+||\vv^n||_{L^1_T(\B^{\frac dp+1}_{p,1})}||a^n||_{\LL^\infty_T(\B^{\frac dp}_{p,1})}
\\&\lesssim \big(X_{p,T}(a^n,\vv^n)\big)^2,
\end{align*}
and
\bbal
&\quad \ ||\nabla(\ln(1+a^n))\nabla \mathbf{v}^n||_{L^1_T(\B^{\frac d2-1}_{2,1})}+||I(a^n)\na a^n||_{L^1_T(\B^{\frac d2-1}_{2,1})}+||\vv^n\cd\na \vv^n||_{L^1_T(\B^{\frac d2-1}_{2,1})}
\\&\lesssim||a^n||_{\LL^2_T(\B^{\frac dp+1}_{p,1})}||\vv^n||_{\LL^2_T(\B^{\frac dp}_{p,1})}+||a^n||_{\LL^\infty_T(\B^{\frac dp}_{p,1})}||\vv^n||_{L^1_T(\B^{\frac dp+1}_{p,1})}+||a^n||_{\LL^\infty_T(\B^{\frac dp-1}_{p,1}\cap \B^{\frac dp}_{p,1})}||a^n||_{L^1_T(\B^{\frac dp+1}_{p,1})}
\\&\quad+||(a^n,\vv^n)||^2_{\LL^2_T(\B^{\frac dp}_{p,1})}+||\vv^n||_{\LL^\infty_T(\B^{\frac dp-1}_{p,1})}||\vv^n||_{L^1_T(\B^{\frac dp+1}_{p,1})}
\\&\lesssim
\big(X_{p,T}(a^n,\vv^n)\big)^2.
\end{align*}
Hence, adding up the above inequalities into \eqref{e4.2}, we get
\bbal
X_{p,T}(a^n,\vv^n)\lesssim X_{p,0}+\big(X_{p,T}(a^n,\vv^n)\big)^2.
\end{align*}
Since $X_{p,0}$ is small enough and $X_{p,T}(a^n,\vv^n)$ depends continuously on the time variable, a standard bootstrap argument will ensure that $T^*_n=+\infty$. Moreover, there holds for $T\in [0,\infty)$,
$$X_{p,T}(a^n,\vv^n)\leq C_0 X_{p,0}, \qquad ||a^n||_{L^\infty_T(L^\infty)}<\frac12.$$
This implies $X_{p,T}(a^n,\vv^n)$ is bounded independent of $n$ for all $T\in[0,\infty)$.

\textbf{Step 3: Existence of the solution.} A classical compactness method as in \cite{Chen 2010} show that we can find a global solution $(a,\vv)\in \mathcal{E}$ satisfying system \eqref{eq1.4} with the initial data $(a_0,\vv_0)$.

\textbf{Step 4: Uniqueness of the solution.} Assume that $(a^1,\mathbf{v}^1)$ and $(a^2,\mathbf{v}^2)$ are two solutions of the system \eqref{eq1.4} with the same initial data $(a_0,\vv_0)$. Setting $\delta a=a^1-a^2$ and $\delta\mathbf{v}=\mathbf{v}^1-\mathbf{v}^2$, we find that $(\delta a,\delta\mathbf{v})$ satisfies
\begin{equation}\label{e4.3}
\begin{cases}
\ \partial_t\delta a-\mu \Delta \delta a+\mathrm{div}\delta \mathbf{v}=f(a^1,\vv^1)-f(a^2,\vv^2),  \\
\ \partial_t\delta \mathbf{v}-\mu\Delta \delta \mathbf{v}+\nabla \delta a=\mathbf{g}(a^1,\vv^1)-\mathbf{g}(a^2,\vv^2),\\
\ (\delta a,\delta \mathbf{v})|_{t=0}=(0,0).
\end{cases}
\end{equation}
Applying the Lemma \ref{lem3.1} to system \eqref{e4.3}, we have for any $T\in[0,\infty)$,
\bal\label{e4.3-1}
{X_{p,T}(\delta a,\delta\vv)}&\lesssim ||\big(f(a^1,\vv^1)-f(a^2,\vv^2),\mathbf{g}(a^1,\vv^1)-\mathbf{g}(a^2,\vv^2)\big)||^\ell_{L^1_T(\B^{\frac d2-1}_{2,1})}\nonumber
\\&\quad+||\big(\na [f(a^1,\vv^1)-f(a^2,\vv^2)],\mathbf{g}(a^1,\vv^1)-\mathbf{g}(a^2,\vv^2)\big)||^h_{L^1_T(\B^{\frac dp-1}_{p,1})}
\\&\lesssim ||f(a^1,\vv^1)-f(a^2,\vv^2)||^\ell_{L^1_T(\B^{\frac d2-1}_{2,1})}+||f(a^1,\vv^1)-f(a^2,\vv^2)||^h_{L^1_T(\B^{\frac dp}_{p,1})}\nonumber
\\&\quad+||\mathbf{g}(a^1,\vv^1)-\mathbf{g}(a^2,\vv^2)||_{L^1_T(\B^{\frac d2-1}_{2,1})}\nonumber.
\end{align}
Note that
\bbal
f(a^1,\vv^1)-f(a^2,\vv^2)=-\D(\delta a \vv^1)-\D(a^2\delta\vv),
\end{align*}
we infer from Lemma \ref{lem2.3} and Lemma \ref{lem2.8} that
\bal\label{e4.3-2}
&\quad ||\D(\delta a\vv^1)+\D(a^2\delta\vv)||^\ell_{L^1_T(\B^{\frac d2-1}_{2,1})}+||\D(\delta a\vv^1)+\D(a^2\delta\vv)||^h_{L^1_T(\B^{\frac dp}_{p,1})}\nonumber
\\&\lesssim ||\vv^1||_{\LL^\infty_T(\B^{\frac dp-1}_{p,1})}||\delta a||_{L^1_T(\B^{\frac dp+1}_{p,1})}+||\vv^1||_{\LL^2_T(\B^{\frac dp}_{p,1})}
||\delta a||_{\LL^2_T(\B^{\frac dp+1}_{p,1}\cap\B^{\frac dp}_{p,1})}\nonumber
\\&\quad+||\vv^1||_{L^1_T(\B^{\frac dp+1}_{p,1})}||\delta a||_{\LL^\infty_T(\B^{\frac dp}_{p,1}\cap\B^{\frac dp-1}_{p,1})}+||\delta\vv||_{\LL^\infty_T(\B^{\frac dp-1}_{p,1})}|| a^2||_{L^1_T(\B^{\frac dp+1}_{p,1})}\\&\quad+||\delta\vv||_{\LL^2_T(\B^{\frac dp}_{p,1})}||a^2||_{\LL^2_T(\B^{\frac dp+1}_{p,1}\cap\B^{\frac dp}_{p,1})}
+||\delta\vv||_{L^1_T(\B^{\frac dp+1}_{p,1})}||a^2||_{\LL^\infty_T(\B^{\frac dp}_{p,1}\cap\B^{\frac dp-1}_{p,1})}\nonumber
\\&\lesssim \ep{X_{p,T}(\delta a,\delta\vv)}.\nonumber
\end{align}
By direct calculation, we also have
\bbal
\mathbf{g}(a^1,\vv^1)-\mathbf{g}(a^2,\vv^2)&=\underbrace{2\mu\nabla[\ln(1+a^1)-\ln(1+a^2)]\nabla \mathbf{v}^1}_{\mathcal{I}_1}+\underbrace{2\mu\nabla(\ln(1+a^2))\nabla \delta\mathbf{v}}_{\mathcal{I}_2}
\\&\quad \underbrace{-I(a^2)\nabla\delta a-[I(a^1)-I(a^2)]\na a^1}_{\mathcal{I}_3}-\underbrace{(\delta\mathbf{v}\cdot\nabla\mathbf{v}^2+\mathbf{v}^1\cdot\nabla\delta\mathbf{v})}_{\mathcal{I}_4}.
\end{align*}
From Lemmas \ref{lem2.3}-\ref{lem2.4}, Remark \ref{re2} and Lemma \ref{lem2.8}, we obtain
\begin{align}\label{e4.4}
||{\mathcal{I}_1}||_{L_T^{1}(\dot{B}_{2,1}^{\frac{d}{2}-1})}
\lesssim&||\ln(1+a^1)-\ln(1+a^2)||_{\LL_T^{\infty}(\dot{B}_{p,1}^{\frac{d}{p}})}||\vv^1||_{L_T^{1}(\dot{B}_{p,1}^{\frac{d}{p}+1})}\nonumber\\
&+||\ln(1+a^1)-\ln(1+a^2)||_{\LL_T^{2}(\dot{B}_{p,1}^{\frac{d}{p}+1})}||\vv^1||_{\LL_T^{2}(\dot{B}_{p,1}^{\frac{d}{p}})}\nonumber\\
\lesssim& || \vv^1||_{L_T^{1}(\dot{B}_{p,1}^{\frac{d}{p}+1})}||\delta a||_{\LL_T^{\infty}(\dot{B}_{p,1}^{\frac{d}{p}})}+|| \vv^1||_{\LL_T^{2}(\dot{B}_{p,1}^{\frac{d}{p}})}||\delta a||_{\LL_T^{2}(\dot{B}_{p,1}^{\frac{d}{p}+1}\cap \dot{B}_{p,1}^{\frac{d}{p}})}\nonumber\\ \lesssim&\ep{X_{p,T}(\delta a,\delta\vv)},
\end{align}
and
\begin{align}\label{e4.5}
||{\mathcal{I}_2}||_{L_T^{1}(\dot{B}_{2,1}^{\frac{d}{2}-1})}
\lesssim&||a^2||_{\LL_T^{\infty}(\dot{B}_{p,1}^{\frac{d}{p}})}||\delta\vv||_{L_T^{1}(\dot{B}_{p,1}^{\frac{d}{p}+1})}+
||a^2||_{\LL_T^{2}(\dot{B}_{p,1}^{\frac{d}{p}+1})}||\delta\vv||_{\LL_T^{2}(\dot{B}_{p,1}^{\frac{d}{p}})}\lesssim\ep{X_{p,T}(\delta a,\delta\vv)}.
\end{align}
To estimate the composition term $\mathcal{I}_3$, it suffices to note that for any sufficiently smooth function $H$, we have
$$H(y)-H(x)=\Big(H'(0)+\int_0^1[H'(x+\tau(y-x))-H'(0)]d\tau\Big)(y-x).$$
Therefore, using $G'(0)=0$ and $\mathcal{I}_3=-\nabla[G(a^1)-G(a^2)]$, we obtain
\begin{align}\label{e4.5.7}
||{\mathcal{I}_3}||_{L_T^{1}(\dot{B}_{2,1}^{\frac{d}{2}-1})}\lesssim &||\nabla\Big(\int_0^1G'(a^1+\tau(a^2-a^1))d\tau\cdot \delta a\Big)||_{L_T^{1}(\dot{B}_{2,1}^{\frac{d}{2}-1})}
\nonumber \\ \lesssim& ||\delta a||_{\LL_T^{2}(\dot{B}_{p,1}^{\frac{d}{p}})}(||a^1||_{\LL_T^{2}(\dot{B}_{p,1}^{\frac{d}{p}})}+||a^2||_{\LL_T^{2}(\dot{B}_{p,1}^{\frac{d}{p}})})
\\\qquad&+||\delta a||_{L_T^{1}(\dot{B}_{p,1}^{\frac{d}{p}+1})}(||a^1||_{\LL_T^{\infty}(\dot{B}_{p,1}^{\frac{d}{p}-1})}+||a^2||_{\LL_T^{\infty}(\dot{B}_{p,1}^{\frac{d}{p}-1})}) \nonumber
\\ \qquad& +||\delta a||_{\LL_T^{\infty}(\dot{B}_{p,1}^{\frac{d}{p}-1})}(||a^1||_{L_T^{1}(\dot{B}_{p,1}^{\frac{d}{p}+1})}+||a^2||_{L_T^{1}(\dot{B}_{p,1}^{\frac{d}{p}+1})})\nonumber
\\ \lesssim&\ep{X_{p,T}(\delta a,\delta\vv)}. \nonumber
\end{align}
For the last term $\mathcal{I}_4$, according to Lemma \ref{lem2.3} and Lemma \ref{lem2.8}, we get
 \begin{align}\label{e4.8}
||{\mathcal{I}_4}||_{L_T^{1}(\dot{B}_{2,1}^{\frac{d}{2}-1})} \lesssim&||\delta\mathbf{v}\cdot\nabla\mathbf{v}^2||_{L_T^{1}(\dot{B}_{2,1}^{\frac{d}{2}-1})}
+||\mathbf{v}^1\cdot\nabla\delta\mathbf{v}||_{L_T^{1}(\dot{B}_{2,1}^{\frac{d}{2}-1})}
\nonumber \\ \lesssim& ||\delta\mathbf{v}||_{\LL_T^{\infty}(\dot{B}_{p,1}^{\frac{d}{p}-1})}||\mathbf{v}^2||_{L_T^{1}(\dot{B}_{p,1}^{\frac{d}{p}+1})}
+||(\mathbf{v}^1,\mathbf{v}^2)||_{\LL_T^{2}(\dot{B}_{p,1}^{\frac{d}{p}})}||\delta \mathbf{v}||_{\LL_T^{2}(\dot{B}_{p,1}^{\frac{d}{p}})}
\\ \qquad& +||\mathbf{v}^1||_{\LL_T^{\infty}(\dot{B}_{p,1}^{\frac{d}{p}-1})}||\delta \mathbf{v}||_{L_T^{1}(\dot{B}_{p,1}^{\frac{d}{p}+1})}
\nonumber
\\ \lesssim &\ep{X_{p,T}(\delta a,\delta\vv)}.\nonumber
\end{align}
Combining the above estimates \eqref{e4.3-1}-\eqref{e4.8} and Lemma \ref{lem3.1}, we can deduce  $${X_{p,t}(\delta a,\delta\vv)}\lesssim \ep{X_{p,t}(\delta a,\delta\vv)},$$
which implies that for all $t\in[0,+\infty)$
$$\mathbf{v}^1(t)=\mathbf{v}^2(t)\quad\mbox{and}\quad a^1(t)=a^2(t).$$
 Therefore, we complete the Proof of Theorem \ref{th1.1}.

\noindent\textbf{Proof of Theorem \ref{th1.2}.} Let $d\geq 3$ and $a_0\in B^{\frac d2}_{2,1}$ and $\vv_0\in B^{\frac d2-1}_{2,1}$. By the Theorem \ref{th1.1} and the construction of system \eqref{eq1.2}, we easily conclude that $(a,\vv)\in L^\infty_T(L^2)$ for any $T\geq 0$ under the assumption of Theorem \ref{th1.1}. Let $(E(t))_{t\geq 0}$ be the semi-group associated with left-hand side of \eqref{eq1.2}, we get after spectral localization
\bbal
U(t)=E(t)U_0+\int^t_0E(t-\tau)F\big(U(\tau)\big)\dd \tau, \quad U(t)=(a(t),\vv(t))\quad\mbox{and}\quad F=(f,\mathbf{g}).
\end{align*}
To simplify the notation, we set
\bbal
&M_\ell(t)=\sup_{0\leq \tau \leq t}(1+\tau)^{\frac d4}\sum_{-1\leq j\leq j_0}(||\Delta_j \vv||_{L^2}+||\Delta_j a||_{L^2}),
\\& M_{h}(t)=\sup_{0\leq \tau \leq t}(1+\tau)^{\frac d4}\sum_{j\geq j_0+1}2^{(\frac d2-1)j}(||\Delta_j \vv||_{L^2}+2^j||\Delta_j a||_{L^2}),
\end{align*}
which implies
\bbal
||a(\tau)||_{B^{\frac d2}_{2,1}}+||\vv(\tau)||_{B^{\frac d2-1}_{2,1}}\leq (1+\tau)^{-\frac d4}M(t), \quad \mathrm{for} \quad \tau\in[0,t].
\end{align*}
where $M(t)=M_\ell(t)+M_h(t)$.\\

To obtain the result of Theorem \ref{th1.2}, we state the estimate of the semi-group operator.

\begin{lemma}\label{le4}
Let $U_0=(a_0,\vv_0)$. Then the operator $E(t)$ satisfies the estimates
\bbal
&\sum_{j\leq j_0}||E(t)\De_{j}U_0||_{L^2}\leq C(1+t)^{-\frac d4}||U_0||_{\B^0_{1,\infty}},
\\&\sum_{j\geq j_0+1}2^{j(\frac d2-1)}||E(t)\De_{j}(\na a_0,\vv_0)||_{L^2}\leq C(1+t)^{-\frac d4}||(\na a_0,\vv_0)||_{\B^{\frac d2-1}_{2,1}}.
\end{align*}
\end{lemma}
{\bf Proof}\quad By Lemma \ref{lem3.1}, it is easy to see that
\bbal
||\mathcal{F}(E(t)\De_jU_0)(\xi)||_{L^2}\leq Ce^{-c_0\mu t2^{2j}}||\De_jU_0||_{L^2}, \quad \mathrm{for} \quad j\leq j_0,
\end{align*}
and
\bbal
||\mathcal{F}(E(t)\De_j(\na a_0,\vv_0))(\xi)||_{L^2}\leq Ce^{-\tilde{c}_0\mu t2^{2j}}||\De_j(\na a_0,\vv_0))||_{L^2}, \quad \mathrm{for} \quad j\geq j_0+1.
\end{align*}
Due to the fact that: for any $\sigma>0$ there exists a constant $C_\sigma$ so that
\bbal
\sup_{t\geq 0}\sum_{k\in \mathbb{Z}}t^{\frac\sigma2}2^{k\sigma}e^{-\frac{c_0}{4} 2^{2k}t}\leq C_{\sigma}.
\end{align*}
Direct calculation shows that
\bbal
t^{\frac d4}\sum_{j\leq j_0}||E(t)\De_jU_0||_{L^2}&\leq C t^{\frac d4}\sum_{j\leq j_0}2^{j\frac d2}e^{-\mu t2^{2j}}||\De_jU_0||_{L^1}
\\&\leq C ||U_0||^\ell_{\B^0_{1,\infty}}\sum_{j\leq j_0}(\sqrt{t}2^j)^{\frac d2}e^{-\mu(\sqrt{t}2^j)^2}\leq C||U_0||_{\B^0_{1,\infty}}
\end{align*}
and
\bbal
\sum_{j\leq j_0}||E(t)\De_jU_0||_{L^2}&\leq C \sum_{j\leq j_0}2^{j\frac d2}e^{-\mu t2^{2j}}||\De_jU_0||_{L^1}
\\&\leq C ||U_0||^\ell_{\B^0_{1,\infty}}\sum_{j\leq j_0}2^{j\frac d2}\leq C||U_0||_{\B^0_{1,\infty}}.
\end{align*}
This implies the first estimate of our result. To obtain the second estimate of this lemma, we find that
\bbal
&t^{\frac d4}\sum_{j\geq j_0+1}2^{j(\frac d2-1)}||E(t)\De_j(\na a_0,\vv_0)||_{L^2}\\
\leq &C t^{\frac d4}\sum_{j\geq j_0+1}2^{j\frac d2}e^{-\mu t2^{2j}}||\De_j(\na a_0,\vv_0)||_{L^2}\\
\leq &C ||(\na a_0,\vv_0)||^h_{\B^0_{2,\infty}}\sum_{j\geq j_0+1}(\sqrt{t}2^j)^{\frac d2}e^{-\mu(\sqrt{t}2^j)^2}\\
\leq &C||(\na a_0,\vv_0)||^h_{\B^0_{2,\infty}}\leq C||(\na a_0,\vv_0)||_{\B^{\frac d2-1}_{2,1}},
\end{align*}
and
\bbal
&\sum_{j\geq j_0+1}2^{j(\frac d2-1)}||E(t)\De_j(\na a_0,\vv_0)||_{L^2}\\
\leq &C \sum_{j\geq j_0+1}2^{j(\frac d2-1)}e^{-\mu t2^{2j}}||\De_j(\na a_0,\vv_0)||_{L^2}\\
\leq &C \sum_{j\geq j_0+1}2^{j(\frac d2-1)}||\De_j(\na a_0,\vv_0)||_{L^2}\\
\leq &C||(\na a_0,\vv_0)||_{\B^{\frac d2-1}_{2,1}}.
\end{align*}
This completes the proof of this lemma.

\noindent Now, we will give the details proof of Theorem \ref{th1.2}. According to Lemma \ref{lem2.6}, Lemma \ref{le4} and the following estimate
\bbal
&\quad \ ||\D(a\vv)||_{L^1}+||\nabla(\ln(1+a))\nabla \mathbf{v}||_{L^1}+||I(a)\na a||_{L^1}+||\mathbf{v}\cdot\nabla\mathbf{v}||_{L^1}
\\&\lesssim ||\vv||_{L^2}||\na a||_{L^2}+||a||_{L^2}||\D \vv||_{L^2}+||\na \vv||_{L^2}||a||_{H^1 \cap B^{\frac d2}_{2,1}}+||a||_{B^{\frac d2}_{2,1}}||a||_{H^1}+||\vv||_{L^2}||\na \vv||_{L^2}
\\&\lesssim (||a||_{B^{\frac d2}_{2,1}}+||\vv||_{B^{\frac d2-1}_{2,1}}) (||\vv||_{\B^{\frac d2+1}_{2,1}}+||a||_{B^{\frac d2}_{2,1}}+||\vv||_{B^{\frac d2-1}_{2,1}}),
\end{align*}
we have for $\tau\in[0,t]$,
\bal\label{eq4.2}
&\quad\int^\tau_0 \sum_{j\leq j_0}||\De_j E(\tau-s)F\big(U(s)\big)||_{L^2}\dd s\nonumber
\\&\leq \int^\tau_0(1+\tau-s)^{-\frac d4}||F\big(U(s)\big)||_{\B^0_{1,\infty}}\dd s
\nonumber\\&\leq \int^\tau_0(1+\tau-s)^{-\frac d4}||F\big(U(s)\big)||_{L^1}\dd s
\\&\leq \int^\tau_0(1+\tau-s)^{-\frac d4}(1+s)^{-\frac d4}M(\tau)\Big((1+s)^{-\frac d4}M(\tau)+||\vv(s)||_{\B^{\frac d2+1}_{2,1}}\Big)\dd s\nonumber
\\&\leq C(1+\tau)^{-\frac d4}M(\tau)X_2(\tau)+C(1+\tau)^{-\frac d2}M^2(\tau)\nonumber,
\end{align}
where
$$X_2(t)=||(a,\vv)||_{\LL^\infty_t(\B^{\frac d2-1}_{2,1})\cap L^1_t(\B^{\frac d2+1}_{2,1})}+||a||_{\LL^\infty_t(\B^{\frac d2}_{2,1})\cap L^1_t(\B^{\frac d2+2}_{2,1})}.$$
Due to the fact that
\bbal
&\qquad||(\na f,\mathbf{g})||_{\B^{\frac d2-1}_{2,1}}
\\&\leq C(||\D(a\vv)||_{\B^{\frac d2}_{2,1}}+||\nabla(\ln(1+a))\nabla \mathbf{v}||_{\B^{\frac d2-1}_{2,1}}+||I(a)\na a||_{\B^{\frac d2-1}_{2,1}}+||\mathbf{v}\cdot\nabla\mathbf{v}||_{\B^{\frac d2-1}_{2,1}})
\\&\leq C(||a||_{\B^{\frac d2}_{2,1}}||\vv||_{\B^{\frac d2+1}_{2,1}}+||\vv||_{\B^{\frac d2-1}_{2,1}}||\vv||_{\B^{\frac d2+1}_{2,1}}+||a||_{\B^{\frac d2-1}_{2,1}}||a||_{\B^{\frac d2+1}_{2,1}}+||a||_{\B^{\frac d2+1}_{2,1}}||\vv||_{\B^{\frac d2}_{2,1}})
\\&\leq C\Big( (||a||_{B^{\frac d2}_{2,1}}+||\vv||_{B^{\frac d2-1}_{2,1}}) ||(a,\vv)||_{\B^{\frac d2+1}_{2,1}}+||a||^{\frac12}_{\B^{\frac d2}_{2,1}}||a||^{\frac12}_{\B^{\frac d2+2}_{2,1}}||\vv||^{\frac12}_{\B^{\frac d2-1}_{2,1}}||\vv||^{\frac12}_{\B^{\frac d2+1}_{2,1}}\Big)
\\&\leq C\Big((||a||_{B^{\frac d2}_{2,1}}+||\vv||_{B^{\frac d2-1}_{2,1}})||(a,\vv)||_{\B^{\frac d2+1}_{2,1}}+||a||_{B^{\frac d2}_{2,1}}||a||_{\B^{\frac d2+2}_{2,1}}\Big),
\end{align*}
we get from Lemma \ref{lem2.6} that for $\tau\in[0,t]$,
\bal\label{eq4.3}
&\quad \ \int^\tau_0||E(\tau-s)F\big(U(s)\big)||^h_{\B^{\frac d2}_{2,1}\times \B^{\frac d2-1}_{2,1}}\dd s\nonumber
\\&\leq \int^\tau_0(1+\tau-s)^{-\frac d4}||(\na f,\mathbf{g})||^h_{\B^{\frac d2-1}_{2,1}}\dd s
\nonumber\\&\leq \int^\tau_0(1+\tau-s)^{-\frac d4}(1+s)^{-\frac d4}M(t)\Big(||(a,\vv)||_{\B^{\frac d2+1}_{2,1}}+||a||_{\B^{\frac d2+2}_{2,1}}\Big)\dd s
\\&\leq C(1+\tau)^{-\frac d4}M(t)X_2(t).\nonumber
\end{align}
By Lemma \ref{le4}, we have
\bal\label{eq4.4}
M(t)&\lesssim M_0+\sup_{0\leq \tau \leq t}(1+\tau)^{\frac d4} \int^\tau_0 \sum_{j\leq j_0}||\De_j E(\tau-s)F\big(U(s)\big)||_{L^2}\dd s\nonumber\\
&~~~~+ \sup_{0\leq \tau \leq t}(1+\tau)^{\frac d4} \int^\tau_0||E(\tau-s)F\big(U(s)\big)||^h_{\B^{\frac d2}_{2,1}\cap\B^{\frac d2-1}_{2,1}}\dd s
\end{align}
Noticing that $X_2(t)\lesssim M_0$ and adding up \eqref{eq4.2}-\eqref{eq4.3} into \eqref{eq4.4}, we conclude that
\bbal
M(t)\lesssim M_0+M(t)M_0+M^2(t).
\end{align*}
If $M_0$ is sufficiently small enough, the standard continuous method show that for all $t\geq 0$,
$$M(t)\leq \widetilde{C}_0M_0,$$
for some positive constant $\widetilde{C}$. Thus, we complete the proof of Theorem \ref{th1.2}.

\vspace*{2em} \noindent\textbf{Acknowledgements. }This work was partially supported by NSFC (No.11801090). 


\begin{thebibliography}{00}\label{ref:ref}\addtolength{\itemsep}{-1.0ex}
\bibitem{Bahouri2011} H. Bahouri, J.Y. Chemin, R. Danchin, Fourier Analysis and Nonlinear Partial Differential Equations, Grundlehren Math. Wiss., vol.343, Springer-Verlag, Berlin, Heidelberg, 2011.
\bibitem{Bresch 2003} D. Bresch, B. Desjardins, C. Lin, On Some Compressible Fluid Models: Korteweg, Lubrication, and Shallow Water Systems, Comm. Part. Diffe. Equ,. 28:3-4 (2003), 843-868,
\bibitem{Chen 2010} Q. Chen, C. Miao, Z. Zhang, Global well-posedness for compressible Navier-Stokes equations with highly oscillating initial velocity. Comm. Pure Appl. Math. 63, (2010) 1173-1224.
\bibitem{C-M-Z} Q. Chen, C. Miao, Z. Zhang, \textit{On the well-posedness for the viscous shallow water equations}, {SIAM J. Math. Anal.}, {\bf40} (2008), 443-474.
\bibitem{C.M.Z} Q. Chen, C. Miao, Z. Zhang, \textit{Well-posedness in critical spaces for the compressible Navier-Stokes equations with density dependent viscosities}, {Rev. Mat. Iberoam.}, {\bf26} (2010), 915-946.
\bibitem{Dunn 1985} J.E. Dunn, J. Serrin, On the thermomechanics of interstitial working. Arch. Ration. Mech. Anal. 88 (1985), 95-133.
\bibitem{Danchin1 2000} R. Danchin, Global existence in critical spaces for compressible Navier-Stokes equations, Invent. Math., 141 (2000),579-614.
\bibitem{Danchin2 2001}  R. Danchin, Local theory in critical spaces for compressible viscous and heat-conductive gases, Comm. Partial Differential Equations, 26 (2001), 1183-1233.
\bibitem{Danchin3 2001} R. Danchin, Global existence in critical spaces for flows of compressible viscous and heat-conductive gases, Arch. Ration. Mech. Anal., 160 (2001), 1-39.
\bibitem{Danchin4 007}R. Danchin, Well-posedness in critical spaces for barotropic viscous fluids with truly not constant density, Comm. Partial Differential Equations, 32 (2007), 1373-1397.
\bibitem{Danchin 2001} R. Danchin, B. Desjardins, Existence of solutions for compressible fluid models of Korteweg type, Ann. Inst. H. Poincar\'{e} Anal. NonLin\'{e}aire, 18 (2001), 97-133.
\bibitem{Danchin 2010} R. Danchin, On the well-posedness of the incompressible density-dependent Euler equations in the $L^p$ framework. Journal of Differential Equations, 248 (2009), 2130-2170.
\bibitem{Danchin5 2016} R. Danchin, L. He, The incompressible limit in $L^p$ type critical spaces, Mathematische Annalen, 366 (2016) 1365-1402.
\bibitem{Danchin6 2017} R. Danchin, P. Muchab, Compressible Navier-Stokes system: Large solutions and incompressible limit, Advances in Mathematics, 320 (2017) 904-925.
\bibitem{Danchin7 2017} R. Danchin, J. Xu, Optimal time-decay estimates for the compressible Navier-Stokes equations in the critical $L^p$ framework, Arch. Rational Mech. Anal., 224 (2017), 53-90.
\bibitem{Haspot 2017} B. Haspot, Global strong solution for the Korteweg system with quantum pressure in dimension $N \geq 2$, Math. Ann. 367: (2017) 667-700.
\bibitem{Haspot 2016} B. Haspot, Existence of global strong solution for Korteweg system with large infinite energy initial data, J. Math. Anal. Appl. 438:1 (2016), 395-443.
\bibitem{Haspot 2011} B. Haspot, Existence of global weak solution for compressible fluid models of Korteweg type, J. Math. Fluid Mech. 13 (2011), 223-249.
\bibitem{Jngel 2010} A. J\"{u}ngel, Global weak solutions to compressible Navier-Stokes equations for quantum fluids, SIAM J. Math. Anal,.42:3 (2010), 1025-1045.
\bibitem{Kawashita} M. Kawashita, On global solution of Cauchy problems for compressible Navier¨CStokes equation, Nonlinear Anal, 48 (2002), 1087-1105.
\bibitem{Li} H.L. Li, T. Zhang, Large time behavior of isentropic compressible Navier¨CStokes system in $\R^3$, Math. Methods Appl. Sci, 34 (2011), 670-682.
\bibitem{Matsumura} A. Matsumura, T. Nishida, The initial value problem for the equation of motion of compressible viscous and heat-conductive fluids, Proc. Japan Acad. Ser, A 55 (1979), 337-342.
\bibitem{Matsumura1980} A. Matsumura, T. Nishida, The initial value problem for the equations of motion of viscous and heatconductive
gases. J. Math. Kyoto Univ, 20 (1980), 67-104.
\bibitem{Okita} M. Okita, Optimal decay rate for strong solutions in critical spaces to the compressible Navier-Stokes equations. Journal of Differential Equations, 257 (2014), 3850-3867.
\bibitem{Ponce}  G. Ponce, Global existence of small solution to a class of nonlinear evolution equations. Nonlinear Anal.TMA, 9 (1985), 339-418.
\bibitem{Wang} Y. Wang, Z. Tan, Global existence and optimal decay rate for the strong solution in $H^2$ to the compressible Navier-Stokes equation, Appl. Math. Lett, 24 (2011), 1778-1784.
\end{thebibliography}
\end{document}